\documentclass[a4,twoside,11pt,leqno,openright]{article}
\usepackage[utf8]{inputenc}	
\usepackage[english]{babel}
\usepackage{latexsym}
\usepackage[T1]{fontenc}
\usepackage{enumerate}
\usepackage{a4wide}
\usepackage{fancyhdr}
\usepackage{amsmath}
\usepackage{amssymb}
\usepackage[all]{xy}
\usepackage{theorem}
\newtheorem{theo}{Theorem}[section]
\newtheorem{defi}{Definition}[section]
\newtheorem{lem}{Lemma}[section]
\newtheorem{rem}{Remark}[section]

\newtheorem{coro}{Corollary}[section]
\newenvironment{preuve}
{\noindent{\textit{Proof~:}}\small}%
{\normalsize\rmfamily\\[.2cm]}
\newenvironment{defini}
{\noindent{\textbf{Definition~}}} %
{\normalsize\rmfamily\\[.2cm]}

\title{Functions with support in a lacunary system of intervals and cyclicity for the semi-group of left translations }
\author{R\'{e}da Choukrallah}
\date{}

\begin{document}

\maketitle
\vspace{1cm}
\begin{abstract}
We prove in the vector-valued spaces $L^2(\mathbb{R}_+,\, X)$ (where $X$ is a finite dimensional Hilbert space) the cyclicity for the semi-group of left translations of some particular functions with support included in a lacunary system of intervals.
\end{abstract}

\section*{Introduction}
Let $L^2(\mathbb{R}_+)$ be the space of measurable functions on $\mathbb{R}_+=[0,\infty)$ such that 
 $$\| f \|_2= \Big ( \displaystyle{\int_0^\infty} |f(x)|^2 dx \Big )^{\frac{1}{2}} < \infty$$
and let $(S_t^*)_{t>0}$ be the semi-group of left translations where  
$$S_t^*f (x)=f(x+t) \, \, x \ge 0, \, \, t>0.$$
If $f \in L^2(\mathbb{R}_+)$, we define $E_f:=span_{L^2(\mathbb{R}_+)} \big (S_t^*f \, : \, t>0 \big )$, and we say that $f$ is cyclic if $E_f=L^2(\mathbb{R}_+)$.\\
We introduce the following definitions for the class of functions we will work with.\\
\\
\begin{defini} Let $X$ be a Hilbert space. A sequence $(a_k)_{k \ge 0}$ of elements in $X$ is said completely relatively compact (c.r.c.) if for any orthogonal projection $P:\, X \rightarrow X$ and $\eta =\{k \ge 0:\, Pa_k \ne 0 \}$, the sequence $\big ( \frac{Pa_k}{\| Pa_k \|} \big )_{k \in \eta}$ is relatively compact in $X$.
\end{defini}
\begin{defini}
A lacunary system of intervals $\Lambda= \displaystyle{\bigcup_{k=1}^{\infty}} [a_k,\, b_k]$ is a countable union of intervals $[a_k, \, b_k]$ verifying the following statements:\\
(1) $a_k <b_k<a_{k+1} \, \, \, \, \forall \, \, k \ge 1$;\\
(2) $\displaystyle{\sup_{k \ge1}} \, \frac{b_k}{a_{k+1}} <1$.
\end{defini}
The main result in this paper is Theorem \ref{rcyc} which states that a function $f \in L^2(\mathbb{R}_+)$ whose support $supp(f)$ is included in a lacunary system of intervals $\Lambda=\displaystyle{\bigcup_{k=1}^{\infty}} [a_k,\, b_k]$ is cyclic in $L^2(\mathbb{R}_+)$ under both conditions
\begin{enumerate}[(1)]
\item $c=\displaystyle{\sup_{k \ge1}} |b_k-a_k|< \infty$;
\item the sequence $( f_k)_{k \ge1}$ where $f_k(t)=f(t+a_k),\, 0 \le t \le c$ is c.r.c. in $L^2(0,\, c)$.
\end{enumerate}
We also generalize this result to the case of vector-valued functions and Theorem \ref{tfinal} gives a criterion of cyclicity for functions in $L^2(\mathbb{R}_+, \, X)$ when statements $(1),\, (2)$ hold.
\section{The scalar case $L^2(\mathbb{R}_+)$}
We first need some well known results in the theory of invariant subspaces and in particular the property of unicellularity for Volterra operator (for completness see \cite{RadRos}).
\begin{defi}
Let $X$ be a Banach space and let $T$ be a bounded operator on $X$.\\
The operator $T$ is said {\em unicellular} if $Lat\,T$ is totally ordered.
\end{defi}
A Volterra-type integral operator on $L^2(0, \, 1)$ is an operator $K$ of the form
$$(Kf)(x)=\int_{0}^{x} k(x,\, y)f(y) \, dy,$$
where $k$ is a square-integrable function on the unit square. The Volterra operator is the particular Volterra-type operator obtained when $k$ is the constant function 1. For any $\alpha \in [0,\, 1]$ let 
$$\mathit{M}_\alpha=\{ f \in L^2(0,\, 1):\, f=0 \,\, \mbox{ a.e. on }\, \, [0,\, \alpha] \}.$$
It is clear that $\{ \mathit{M}_\alpha:\, \alpha \in [0,\, 1] \} \subset Lat \, K$ for any Volterra-type integral operator.\\
The following theorem plays an important role in the construction of our result (see \cite[p.\,68]{RadRos}). 
\begin{theo}\label{Volt}
Volterra operator is unicellular and 
$$Lat \, V=\{ \mathit{M}_\alpha:\, \alpha \in [0,\, 1 ]\}.$$
\end{theo}
To connect the semi-group of translations $(S_t^*)_{t>0}$ with the Volterra operator, we recall the following  result from \cite{Ni3}.
\begin{lem}
If $M$ is a closed subspace in $L^2(0, \, a)$ and $S_t$ is the semi-group of right translations then,
$$S_tM \subset M, \, \, \, \, \forall \, \, t>0 \Leftrightarrow VM \subset M.$$
\end{lem}
The description of invariant subspaces for the Volterra operator, its unicellularity and the connection with the semi-group of right translations leads to the next theorem which is a direct consequence of the previous facts.
\begin{theo}\label{agdo}
$f\in L^2(0,\, a),\, \, a \in supp(f) \,$ $\,\Rightarrow$ $\, E_f=L^2(0,\, a)$.
\end{theo}
\begin{rem}
In what follows, it is important to see that we embedded $L^2(0,\, a)$ in $L^2(0,\, \infty)$ by extending any function of $L^2(0,\, a)$ in $(a, \, \infty)$ by 0.
\end{rem}
We can now start the construction of our proof. 
\begin{lem}\label{modif1}
Let $\Lambda= \displaystyle{\bigcup_{k=1}^{\infty}} [a_k,\, b_k]$ be a lacunary system of intervals on $\mathbb{R}_+$.\\
The function $\phi(x)=card \big \{ (a_j-a_k,\, b_j-a_k):\, x \in (a_j-a_k,\, b_j-a_k) \big \}$ is bounded on $\mathbb{R}_+$.
\end{lem}
\begin{preuve}
It is a slight adaptation from the discret case to the continous one (see \cite{DSS}).
Since $\Lambda$ is a lacunary system of intervals, there exists $d>1$ such that $\displaystyle{\inf_{k \ge1}} \, \frac{a_{k+1}}{b_k} \ge d$. Let $M>0$ be such that $d^{M-1}(d-1) \ge 1$ and let $x \in \mathbb{R}_+$. Then there exists a finite set $J$ such that $\forall \,j \in J,\,\, \, a_j-a_k\le x \le b_j-a_k$.\\
Denoting by $j_0$ the smallest index in $J$ such that $b_{j_0} >x$, we prove that if $j \in J$ then $j_0 \le j < j_0+M$.\\
The first inequality is trivial, while for the second one we argue by contradiction. Suppose that $j \ge j_0+M$. This implies
$$
a_j-a_k \ge a_j-b_{j-1} \ge (d-1)b_{j-1} \ge (d-1)d^{M-1}b_{j_0}>x.
$$
If $a_j-a_k>x$ then $x \not \in (a_j-a_k,b_j-a_k)$ which contradicts that $j \in J$. Consequently, $\phi(x) \le M$ for every $x \in \mathbb{R}_+$ and thus $\phi$ is bounded.
\hfill{$\Box$} 
\end{preuve}
\\[3mm]
Let $f \in L^2(\mathbb{R}_+)$ be such that $supp(f) \subset \Lambda$. Then $f$ can be written \\
\begin{equation*}
f=\displaystyle{\sum_{k \ge 1}} g_k \,,\, \, \, \mbox{ with }supp(g_k) \subset [a_k,\, b_k],
\end{equation*}
or, equivalently,
\begin{equation*}
f(x)=\displaystyle{\sum_{k \ge 1}} f_k(x-a_k), \, \, \, \mbox{ where }f_k \in L^2(0,\, b_k-a_k) \subset L^2(\mathbb{R}_+).\\
\end{equation*}
\begin{theo}\label{rcyc}
Let $f \in L^2(\mathbb{R}_+)$ and $supp(f) \subset \Lambda$, where $\Lambda= \displaystyle{\bigcup_{k=1}^{\infty}} [a_k,\, b_k]$ is a lacunary system of intervals. Suppose that
\begin{enumerate}[i)]
\item $c=\displaystyle{\sup_{k \ge1}} |b_k-a_k|< \infty$;
\item $\big( \frac{f_k}{\| f_k \|} \big)_{k \ge1}$ is relatively compact.
\end{enumerate}
Then $f$ is cyclic in $L^2(\mathbb{R}_+)$.
\end{theo}
\begin{preuve}
For any $k \ge 1$,
$$\displaystyle{\frac{S_{a_k}^*f}{\| f_k \|}}=\frac{f_k}{\| f_k \|}+\displaystyle{\sum_{j > k}} \frac{f_j}{\| f_k \|}.$$
Let $r_k=\displaystyle{\sum_{j > k}} \frac{f_j}{\| f_k \|}$. We use in the first step the same approach developped by R. G. Douglas, H.S. Shapiro, and A.L. Shields (see \cite{DSS}). We consider neighborhoods of the form 
$$\mathit{V}=\big \{ h \in L^2 (\mathbb{R}_+ )\, :\, |(h,h_i)|< 1,\, 1\le i \le n\big
\},$$
where $h_i, \, (i=1 \dots n)$ are fixed functions in $ L^2(\mathbb{R}_+)$.\\
Suppose that none of the $r_k$ is in $\mathit{V}$. Then $1 \le \displaystyle \max_{1 \le i \le n}|(r_k,\, h_i)| \, \, \, \, \forall \, \, k \ge 1.$\\
We obtain that \\
$$|(r_k,\, h_i)| \le \frac{1}{\| f_k \|}\displaystyle{\int_{a_j \le x \le b_j}} |\displaystyle{\sum_{j>k}} f_j(x-a_j+a_k)h_i(x) |dx
\le \frac{1}{\| f_k \|}\Big (\displaystyle{\sum_{j>k}}\| f_j \|^2 \Big)^{\frac{1}{2}} \Big (\displaystyle{\sum_{j>k}}\displaystyle{\int_{a_j-a_k} ^{b_j-a_k}} |h_i(x)|^2 dx \Big)^{\frac{1}{2}}.$$\\
Therefore,\\
$$\frac{\| f_k\|^2}{\displaystyle{\sum_{j>k}}\| f_j \|^2} \le \max_{1 \le i \le n}\displaystyle{\sum_{j>k}}\displaystyle{\int_{a_j-a_k} ^{b_j-a_k}} |h_i(x)|^2 dx \le \displaystyle{\sum_{i=1}^{n}}\displaystyle{\sum_{j>k}}\displaystyle{\int_{a_j-a_k} ^{b_j-a_k}} |h_i(x)|^2 dx$$
and so, by summing on $k$,\\
$$\displaystyle{\sum_{k\ge1}} \frac{\| f_k\|^2}{\displaystyle{\sum_{j>k}}\| f_j \|^2} \le \displaystyle{\sum_{i=1}^{n}}\displaystyle{\sum_{k\ge1}} \displaystyle{\sum_{j>k}}\displaystyle{\int_{a_j-a_k} ^{b_j-a_k}} |h_i(x)|^2 dx.$$
The left-hand side of the above inequality diverges. It suffices to prove that the right-hand side is bounded to have a contradiction and to conclude. The second step is to show that  $$\displaystyle{\sum_{k\ge1}}\displaystyle{\sum_{j>k}}\displaystyle{\int_{a_j-a_k} ^{b_j-a_k}} |h_i(x)|^2 dx < \infty.$$
Suppose $\phi(x)=\displaystyle{\sum_{k\ge1}}\displaystyle{\sum_{j>k}}\chi_{(a_j-a_k,b_j-a_k)}(x)$, and note that the function $\phi$ can also be written in the following form, $$\phi(x)=card \big \{ (a_j-a_k,b_j-a_k):x \in (a_j-a_k,\, b_j-a_k) \big \}.$$ 
According to Lemma \ref{modif1}, $\phi$ is bounded, so there exists $M>0$ such that $\phi(x) \le M$ for every $x \in \mathbb{R}_+$. It follows that\\
\begin{eqnarray*}
\displaystyle{\sum_{k\ge1}}\displaystyle{\sum_{j>k}}\displaystyle{\int_{a_j-a_k} ^{b_j-a_k}} |h_i(x)|^2 dx & = &\displaystyle{\sum_{k\ge1}}\displaystyle{\sum_{j>k}}\displaystyle{\int_{\mathbb{R}_+}} \chi_{(a_j-a_k,b_j-a_k)}(x) |h_i(x)|^2 dx\\
&= &\displaystyle{\int_{\mathbb{R}_+}} \displaystyle{\sum_{k\ge1}}\displaystyle{\sum_{j>k}}\chi_{(a_j-a_k,\,b_j-a_k)}(x) |h_i(x)|^2 dx\\
&  = &\displaystyle{\int_{\mathbb{R}_+}} \phi(x) |h_i(x)|^2 dx \le M \displaystyle{\int_{\mathbb{R}_+}} |h_i(x)|^2 dx < \infty. \\
\end{eqnarray*}
Then each weak neighborhood of $0$ contains at least one of the $(r_{k})_{k \ge1}$. From
$$S_{a_{k_j}}^*f=\frac{f_{k_j}}{\|f_{k_j}\|}+r_{k_j},$$ 
and since $( f_k / \| f_k \|)_{k \ge1}$ is relatively compact, there is a subsequence $(k_{j_i})_{i \ge 1}$, such that $(f_{k_{j_i}}/\|f_{k_{j_i}}\|)_{i \ge 1}$ converges weakly to a limit $g \ne 0$ and therefore $g \in E_f.$\\
Setting $b=sup \, \, supp(g)$, condition $(i)$ gives that $0<b<\infty$ and $g \in L^2(0,\, b)$.\\
Note that
$E_g \subset E_f,$ 
because $S_t^* E_f \subset E_f \, \, \, \, \forall \, \, t>0$.
We apply Theorem \ref{agdo} and,
$$E_g=L^2(0,\,b) \subset  E_f.$$
Reasoning as before for $S_{a_k-T}^*$ where $T>0$ is arbitrary and $k$ is such that $a_k>T$, we obtain  $$S_{a_k-T}^*f(x)=\frac{f_k(x+T)}{\|f_k\|}+r_k(x).$$
We extract again a subsequence, $(f_{k_{j_i}}/ \|f_{k_{j_i})\| } )_{i \ge 1}$ whose limit $g$ is a nonzero function such that $g \in E_f$ and $supp(g) \subset [T,\, c+T]$. Since $b=sup \, \, supp(g)>T$ and by Theorem \ref{agdo} we obtain
 $$L^2(0,\, T) \subset E_f \, \, \, \, \forall \, \, T>0.$$ 
Finally, $L^2(\mathbb{R}_+) \subset E_f$.
\hfill{$\Box$}
\end{preuve}
\section{The vector-valued case $L^2(\mathbb{R}_+, \, \mathbb{C}^n)$}
In this section, we generalize the previous result to $L^2(\mathbb{R}_+,\, X)$ where $X$ is a finite dimensional Hilbert space. Let $c>0$ be a constant and $(f_k)_{k \ge 1}$  be a c.r.c. sequence in $L^2((0,\, c),\, X)$. Then the sequences $\big ( \frac{Pf_k}{\|Pf_k \|} \big )_{k \in \eta}$ are relatively compact in $L^2((0,\, c),\, PX)$ for any orthogonal projection $P:\, L^2(0,\, c;\, X) \rightarrow L^2(0,\, c;\, X)$, where $\eta=\{k \ge 1:\, Pf_k \ne 0\}$.\\[3mm]
The next lemma is the first step for the generalization to vector-valued functions.
\begin{lem}\label{rlem}
Let $f \in L^2(\mathbb{R}_+,\, X)$ and $supp(f) \subset \Lambda= \displaystyle{\bigcup_{k=1}^{\infty}} [a_k,\, b_k]$ where $\Lambda$ is a lacunary system of intervals. For any $k \ge0$, we consider $f_k(x)=f(x+a_k),$ with $supp(f_k) \subset [0,\, b_k-a_k]$ and suppose that,
\begin{enumerate}[i)]
\item $\displaystyle{\sup_{k \ge1}} |b_k-a_k| \le c < \infty$.
\item $\big( \frac{f_k}{\| f_k \|} \big)_{k \ge 1}$ is relatively compact in $L^2(0,\, c; \, X)$.
\end{enumerate}
Then, there exists a nonzero function $g \in L^2(0,\, c;\, X)$ such that $S_tg \in E_{f}$ pour tout $t\ge 0$.
\end{lem}
\begin{preuve}
It is almost the same as the proof of the first part of Theorem  \ref{rcyc} with some minor technical modifications and we change the products $f_j(x-a_j+a_k)h_i(x)$ by the scalar products $(f_j(x-a_j+a_k),\, h_i(x))_X$.
\hfill{$\Box$}
\end{preuve}
\begin{rem}
In fact, we can prove more, however it is not necessary for what follows. The function $g$ verifies 
$$S_tg \in E_{S_u^*f} \, \, \, \, \forall \, \, t,\, u \ge 0.$$
\end{rem}
Let $c>0$, and let $M$ be a subspace of $L^2(0,\, c,\, X)$. Recall that any space $L^2(a,\,b;\,X)$ can be seen as a subspace of $L^2(\mathbb{R}_+,\,X)$. We take
$$L^2\otimes M=span_{L^2(\mathbb{R}_+,\,X)}(S_tM:\, t \ge 0)=clos_{L^2(\mathbb{R}_+,\,X)}\Big\{\displaystyle{\sum_{i}}\varphi_i*m_i:\, m_i \in M,\, supp(\varphi_i)\mbox{ compact}\Big\},$$
where $\varphi_i$ are functions in $L^2$ or $L^1$ (there is no influence on the result and the sums are finite). With this notation, we have the following.
\begin{coro}\label{corovv}
Let $f$ as in Lemma \ref{rlem}. \\
There exists a closed subspace $M=M(f) \subset L^2(0,\, c;\, X)$ such that
\begin{enumerate}[(a)]
 \item $M(f) \ne \{0 \}$.
 \item $L^2\otimes M(f) \subset E_f$.
 \item If $h \in L^2(0,\, c;\, X)$ and $L^2\otimes h \subset E_f$, then $h \in M(f)$, (``$M(f)$ is ``maximal'').
\end{enumerate}
\end{coro}
\begin{preuve}
Lemma \ref{rlem} gives the existence of one closed subspace $M$ satisfying $(a)$ and $(b)$.  According to the definition of the ''product`` $L^2\otimes M$, it is clear that if $L^2\otimes M_\alpha \subset E_f\, (\alpha \in A)$, then $L^2\otimes M \subset E_f$ where $M=span_{L^2(0,\, c;\, X)}(M_\alpha:\, \alpha \in A)$, which proves property $(c)$. 
\hfill{$\Box$}
\end{preuve}
The following theorem contains other properties of the subspace $L^2\otimes M(f) \subset E_f$ in Corollary \ref{corovv} which are connected with the possible cyclicity of the function $f$.
\begin{theo} \label{teovv}
Let $f$ as in Lemma \ref{rlem} such that the sequence $(f_k)_{k \ge 1}$ is c.r.c. in $L^2(0,\, c;\, X)$ and let $F=L^2\otimes M(f)$ be the subspace in $E_f$ defined in Corollary \ref{corovv}. Then
\begin{enumerate}[(i)]
\item $S_tF \subset F \subset E_f$ for any $t,\, t \ge 0$.
\item $f=g+h$ where $g\in F$ and $supp(h)$ is compact.
\item Let $E_{M(f)}=span_{L^2(\mathbb{R}_+,\,X)}(S_t^*M(f):\, t \ge 0)$ be the $S_t^*$-invariant subspace generated by $M(f)$. Then, $E_{M(f)} \subset L^2(0,\,c; \,X)$ and 
$$E_f=clos_{L^2(\mathbb{R}_+,\,X)}(L^2\otimes M(f)+E_{M(f)}).$$
The subspaces $L^2\otimes M(f)$ and $E_{M(f)}$ are of the form 
$$L^2\otimes M(f)=\mathcal{F}(\Theta_1H^2(\mathbb{C}_+,\,X_1)),\, E_{M(f)}=L^2(\mathbb{R}_+,\,X)\ominus \mathcal{F}(\Theta_2 H^2(\mathbb{C}_+,\, X_2)),$$
where $X_i \subset X(i=1,\, 2)$, $\mathcal{F}$ is the Fourier transform, $\Theta_i$ are left inner matricial functions, $\Theta_i(z):\, X_i \rightarrow X$ and $\Theta_2$ is a left divider of $e^{icz}\, I_X.$
\end{enumerate}
\end{theo}
\begin{preuve}
Assertion $(i)$ is obvious by the definition of $F$. In order to prove $(ii)$, observe that the sum 
$$g=\displaystyle{\sum_{k \ge 1}} S_{a_k}\varphi_k,\, \ \  \varphi_k \in M,$$
with convergence in $L^2(\mathbb{R_+},\,X)$ is in $L^2\otimes M$. In particular, it is the case for $\varphi_k=P_Mf_k$ where the $f_k$ are the functions of the lacunary decomposition of $f$, $P_M$ is the orthogonal projection on $M$ in $L^2(0,\, c;\, X)$. Then, the function
$$h=f-g=\displaystyle{\sum_{k \ge 1}} S_{a_k}P_{M^\bot}f_k$$
is in $E_f$, where $M^\bot=L^2(0,\, c;\, X)\ominus M$. It is clear that all limits of the sequence $(P_{M^\bot}f_k/\|P_{M^\bot}f_k \|)_{k \ge 1}$ are in $M^\bot$. Therefore, the support $supp(h)$ is compact or otherwise using the c.r.c. property of $(f_k)_{k \ge 1}$, and by applying again Lemma \ref{rlem} to the function $h$ we obtain a nonzero function $\varphi \in E_f\cap M^\bot$ such that $S_t\varphi \in E_f$ for any $t \ge 0$; which gives a contradiction with the definition of $M=M(f)$.\\
To prove $(iii)$, we note that $L^2\otimes M(f) \subset E_f$ and $E_{M(f)} \subset E_f \cap L^2(0,\, c;\, X)$, which lead to
$$clos_{L^2(\mathbb{R}_+,\,X)}(L^2\otimes M(f)+E_{M(f)}) \subset E_f.$$
In order to establish the converse inclusion, we consider a finite combination $g=\sum_k S_{a_k}\varphi_k, \, \varphi_k \in M$. For any $t \ge 0$, we have
$$S_t^*g=\displaystyle{\sum_{a_k<t}} S_t^*S_{a_k}\varphi_k+\displaystyle{\sum_{a_k\ge t}} S_t^*S_{a_k}\varphi_k=\displaystyle{\sum_{a_k<t}} S_{t-a_k}^*\varphi_k+\displaystyle{\sum_{a_k\ge t}} S_{a_k-t}\varphi_k,$$
where $\displaystyle{\sum_{a_k<t}} S_{t-a_k}^*\varphi_k$ is in $E_{M(f)}$ and $\displaystyle{\sum_{a_k\ge t}} S_{a_k-t}\varphi_k$ in $L^2\otimes M(f)$. This shows that the $S_t^*$-invariant subspace generated by $L^2\otimes M(f)$ and $E_{M(f)}$ is in $clos_{L^2(\mathbb{R}_+,\,X)}(L^2\otimes M(f)+E_{M(f)})$. But this subspace contains also $E_f\,$, so $E_f \subset clos_{L^2(\mathbb{R}_+,\,X)}(L^2\otimes M(f)+E_{M(f)})$.\\
The representations for $L^2\otimes M(f)$ and $E_{M(f)}$ come from Lax-Halmos Theorem and the fact that the inclusion  
$$E_{M(f)}=L^2(\mathbb{R}_+,\,X)\ominus \mathcal{F}(\Theta_2 H^2(X_2))\subset L^2(0,\, c;\, X)=L^2(\mathbb{R}_+,\,X)\ominus \mathcal{F}(e^{icz}H^2(X))$$
is equivalent to the mentioned division in the statement (see, e.g., \cite[p.\,19]{Ni1}).
\hfill{$\Box$}
\end{preuve}
We give the Fourier transforms of the previous objects,
$$F=\mathcal{F}^{-1}(f) \in H^2(\mathbb{C}_+,\, X),\, \mathcal{E}_F=\mathcal{F}^{-1}(E_f)=span_{H^2(\mathbb{C}_+,\, X)}(P_+e^{-itz}F:\, t \ge 0),$$
$$\mathcal{M}(f)=\mathcal{F}^{-1}(M(f)) \subset \mathcal{F}^{-1}(L^2(0,\, c;\, X))=K_\theta=H^2(\mathbb{C}_+,\, X)\ominus \theta H^2(\mathbb{C}_+,\, X),$$
where $\theta=e^{icz}$,
$$H_+^2\otimes \mathcal{M}(f):=\mathcal{F}^{-1}(L^2\otimes M(f))=span_{H^2(\mathbb{C}_+,\, X)}(e^{itz}\mathcal{M}(f):\, t \ge 0)=\Theta_1 H^2(\mathbb{C}_+,\, X_1),$$
$$\mathcal{E}_{\mathcal{M}(f)}=\mathcal{F}^{-1}(E_{M(f)})=K_{\Theta_2}:=H^2(\mathbb{C}_+,\, X)\ominus \Theta_2 H^2(\mathbb{C}_+,\, X_2).$$
\begin{coro}\label{corovv2}
Under the hypotheses of Theorem \ref{teovv} and with the previous notation, we have  
$$\mathcal{E}_F=clos_{H^2(\mathbb{C}_+,\, X)}(\Theta_1 H^2(\mathbb{C}_+,\, X_1)+K_{\Theta_2}),$$
and $dim\, X_1=Rank(\mathcal{M}(f))$, where
$$Rank(\mathcal{M}(f))=max(dim \{\Phi(z):\, \Phi \in \mathcal{M}(f)\}:\, Im(z)>0).$$
\end{coro}
\begin{preuve}
The first formula is the Fourier transform of the one in part $(iii)$ of Theorem \ref{teovv}. The second comes from 
$$dim\, X_1=\displaystyle{\max_z}(dim\{\Phi(z):\, \Phi \in \Theta_1 H^2(\mathbb{C}_+,\, X_1)\},$$
and from the definition of $\Theta_1,\, span_{H^2(\mathbb{C}_+,\, X)}(e^{itz}\mathcal{M}(f):\, t \ge 0)=\Theta_1 H^2(\mathbb{C}_+,\, X_1).$
\hfill{$\Box$}
\end{preuve}
\begin{lem}\label{lemvv}
Let $f\in L^2(\mathbb{R}_+,\,X)$ be as in the statement of Theorem \ref{teovv}, and suppose that $dim\, X<\infty$. If $Rank(\mathcal{M}(f))=dim\,X$, where $\mathcal{M}(f)=\mathcal{F}^{-1}(M(f))$, then $E_f=L^2(\mathbb{R}_+,\,X)$ ($f$ is cyclic).
\end{lem}
\begin{preuve}
If $dim\, X_1=dim\, X$ then $X_1=X$, and 
$$(det\, \Theta_1)H^2(\mathbb{C}_+,\, X) \subset  \Theta_1 H^2(\mathbb{C}_+,\, X)\subset \mathcal{E}_F.$$
We have that $P_+e^{-itz}\mathcal{E}_F \subset \mathcal{E}_F,\, t \ge 0$, $\theta=det(\Theta_1)$ is a scalar inner function and the $P_+e^{-itz}$-invariant subspace generated by $\theta H^2(\mathbb{C}_+,\, X)$ coincides with $H^2(\mathbb{C}_+,\, X)$. Indeed, if $\phi_n$ are the Fejer polynomials of $\theta$, then $$\displaystyle{\lim_{n \rightarrow \infty}}\|P_+\overline{\phi}_n\theta f-f\|_2=\displaystyle{\lim_{n \rightarrow \infty}}\|\phi_n(S^*)\theta f-f\|_2=0,$$ 
for every $f \in H^2(\mathbb{C}_+,\, X)$, and therefore, $\mathcal{E}_F=H^2(\mathbb{C}_+,\, X)$.
\hfill{$\Box$}
\end{preuve}
\begin{theo} \label{tfinal}
Let $f \in L^2(\mathbb{R}_+,\,X)$ satisfy the conditions of Theorem \ref{teovv}, and suppose that $dim\, X<\infty$. Then, $f$ is cyclic for the semi-group $(S_t^*)_{t \ge 0}$ in $L^2(\mathbb{R}_+,\,X)$ if and only if
$$Rank(\mathcal{M}(f))=dim\, X,$$
where $\mathcal{M}(f)=\mathcal{F}^{-1}(M(f))$.
\end{theo}
\begin{preuve}
The sufficiency part was proved in Lemma \ref{lemvv}. Conversely, if $f$ is a cyclic function then
$$\mathcal{E}_F=H^2(\mathbb{C}_+,\, X)=clos_{H^2(\mathbb{C}_+,\, X)}(\Theta_1 H^2(\mathbb{C}_+,\, X_1)+K_{\Theta_2}),$$
by Corollary \ref{corovv2}. Taking $\theta=e^{icz}$ and since $K_{\Theta_2} \subset K_\theta$, we obtain
$$H^2(\mathbb{C}_+,\, X)=P_+\overline{\theta}(H^2(\mathbb{C}_+,\, X))=clos_{H^2(\mathbb{C}_+,\, X)}(P_+\overline{\theta}\Theta_1 H^2(\mathbb{C}_+,\, X_1)).$$
This implies that the subspace $\Theta_1H^2(\mathbb{C}_+,\, X_1)$ is $(P_+e^{-itz})_{t \ge 0}$-cyclic in $H^2(\mathbb{C}_+,\, X)$.\\
Suppose now that $Rank(\mathcal{M}(f))<dim\, X$. Then, using again Corollary \ref{corovv2}, $dim\, X_1 < dim\, X$. According to the Lemma of Complementation (see \cite{Ni1}), there exists a (bilateral) inner function, $\Theta(z): \, X \longrightarrow X$ such that $\Theta \vert_{X_1}=\Theta_1$. Taking into account that the multiplication by a left inner function is an isometric application, we have 
$$\Theta_1 H^2(\mathbb{C}_+,\, X_1)\perp \Theta^\prime H^2(\mathbb{C}_+,\, X^\prime),$$
where $X^\prime=X\ominus X_1 \ne \{0\}$ and $\Theta^\prime=\Theta \vert_{X^\prime}$.\\
Therefore, $\Theta_1H^2(\mathbb{C}_+,\, X_1)\perp e^{itz}\Theta^\prime H^2(\mathbb{C}_+,\, X^\prime)$ for any $t \ge 0$, which leads to a contradiction with the cyclicity of $\Theta_1 H^2(\mathbb{C}_+,\, X_1)$. Finally $Rank(\mathcal{M}(f))=dim\, X$.
\hfill{$\Box$}
\end{preuve}
\begin{rem}
The $S_t^*$-invariant subspace  $E_{M(f)}$ which is included in $L^2(0,\,c;\, X)$ has the Fourier transform such that  
$$E_{M(f)}=L^2(0,\,c;\, X)\ominus \mathcal{F}(\Theta H^2(X)),$$
where $\Theta$ is an inner function with the spectrum $\{0\}$. This leads after the change of variable $z=(x-i)/(x+i)$ and according to the Theorem of factorization (see \cite{GiSh}) to have,
$$\Theta(z)=\int_a^{\stackrel{\curvearrowleft}{b}} exp \bigg(\frac{z+1}{z-1}\bigg)\, d\mu(t),$$
where $\mu$ is a positive matricial measure.
\end{rem}
We end with a simple and instructive example.\\
\\
\textbf{Example:} Let $X=\mathbb{C}^2$, $a_k=2^k,\, b_k=2^k+2,\, f=(\varphi_1,\, \varphi_2)$ where $\varphi_1=\sum_kg_k\chi_{\Delta_k}$, with $g_k$ nonzero continuous functions on $\Delta_k=[2^k+1,\, 2^k+2]$ such that $\sum_k\|g_k\chi_{\Delta_k}\|<\infty$, and $\varphi_2=S_1^*\varphi_1$. With convenient functions $g_k$, we can have the space $M(f)$ of infinite dimension, but with $Rank(\mathcal{F}M(f))=1$. The function $f$, of course, is not cyclic.
\hfill{$\blacksquare$}


\begin{thebibliography}{99}
\bibitem{DSS}
R. G. \textsc{Douglas}, H. S. \textsc{Shapiro}, A. L. \textsc{Shields},
\emph{Cyclic vectors and invariant subspaces for the backward shift operator}, Ann. Inst. Fourier (Grenoble), 20 (1970), fasc.\,1, 37--76.
\bibitem{GiSh}
Yu. P. \textsc{Ginzburg}, L. V. \textsc{Shevchuk},
\emph{On the Potapov theory of multiplicative representations. Matrix and operator valued functions},
Oper. Theory Adv. Appl., Birkh\''auser, Basel, 1994, 72,  28--47.
\bibitem{Ni1}
N. K. \textsc{Nikolskii}, 
\emph{Treatise on the shift operator},
Springer-Verlag, Berlin, 1986. 
\bibitem{Ni3}
N. K. \textsc{Nikolskii}, 
\emph{Invariant subspaces in operator theory and function theory},
Mathematical analysis, Vol. 12 (Russian),  pp. 199--412, 468. Akad. Nauk SSSR Vsesojuz. Inst. Nau\v cn. i Tehn. Informacii, Moscow, 1974.; Engl. transl.: J. Soviet. Math., 5 (1976), 129--249.
\bibitem{RadRos}
H. \textsc{Radjavi}, P. \textsc{Rosenthal} 
\emph{Invariant Subspaces},
Springer Verlag, 1973.
\end{thebibliography}
\end{document}